# Hirsch-type equations and bundles

Leo Egghe

Hasselt University, Belgium

E-mail: leo.egghe@uhasselt.be

Abstract.

We define Hirsch-type equations and bundles being common generalizations of the defining equations of e.g. Hirsch-bundles, g-bundles and Kosmulski-bundles. In this way, common properties of alle these bundles can be proved. The main result proves basic inequalities for these bundles. They form the basis for convergence results as well as for criteria for these bundles to be impact bundles.

## 1. Introduction

We start with the most famous definition in informetrics in the latest two decades, namely the $h$-index (or Hirsch-index). It is defined as follows.

Let $f: \mathbb{R}^+ \to \mathbb{R}^+$ be a decreasing function, representing a rank-frequency function in informetrics of e.g., an author. Then $f(x)$ is the number of items (e.g., citations) in the source (e.g., an article) on rank $x$. Then $h$ is the value $x$, if it exists, such that

$$f(x) = x \quad (1)$$

Note that $h$ is the abscissa of the intersection of the graph of $f$ and the straight line $y = x$. We refer to Hirsch (2005) for the original definition in a discrete setting, see also (Rousseau, Egghe and Guns, 2018). The $h$-index is illustrated in Fig.1

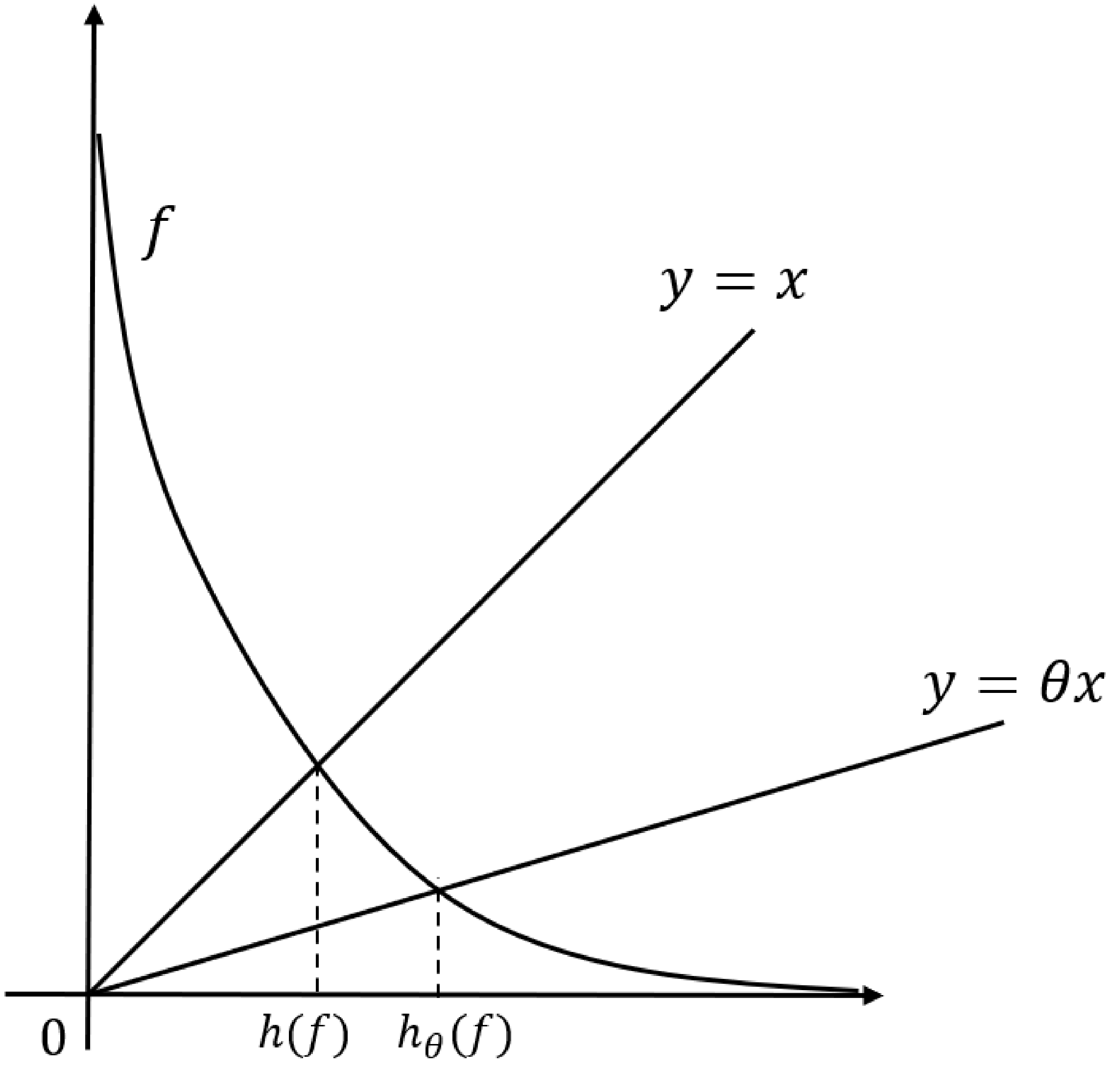

Fig. 1. $h$-indices for $\theta = 1$ and for a general $\theta$.

Since the function $f$ is decreasing, $h$ is unique if it exists. This means that the intersection of the graph of $f$ and the straight line $y = x$ is not empty. The index $h$ can be seen as a measure of impact of the function $f$ (more on this in the sequel).

Definition (1) is remarkable since it is a functional relation derived from $f$. Of course, as remarked by van Eck and Waltman (2008), see also (Egghe & Rousseau, 20XX) we do not have to intersect $f$ with $y = x$ but any straight line $y = \theta x$, $\theta > 0$, can be used. Even more general functions of $\theta$ and $x$ can be used (see further on). Now we have the defining relationship of the $\theta$- dependent $h$ –index $h_\theta(f)$ of $f$ :

$$f(x) = \theta x \tag{2}$$

where again $x = h_\theta(f)$ is unique if it exists (see also Fig. 1). Note that $h_1(f) = h$. The above definition hence yields a new informetric function, called the Hirsch function:

$$h(f): \theta \rightarrow h_\theta(f)$$

which is an important tool in informetrics. For this, we refer the reader to Egghe (2023), a study of its properties. Note that this function always exists on $\mathbb{R}^+$ if $f$ is continuous on an interval $[0,S], S > 0$, with f(S) = 0, or on $\mathbb{R}^+$, with $\lim_{x\to+\infty} f(x) = 0$. If $f(S) > 0$, then $h_\theta(f)$ exists for $\theta \in \left[\frac{f(S)}{S}, +\infty\right]$, and similarly for other subintervals of $\mathbb{R}^+$. In all cases $h_\theta(f)$ is unique. It is well-known that, although $h(f)$ is a very basic notion in informetrics, the (generalized) $h$-index has some pitfalls in the impact evaluation of the decreasing function $f$, see (Egghe, 2006a,b; Rousseau et al., 2018), the most important one being the fact that $h$ is insensitive for the number of items in the most productive sources. Therefore, in (Egghe, 2006a,b) this author introduced the $g$-index, defined as a variant of (1). Here we immediately introduce the variant of equation (2).

The $\theta$-dependent $g$-index of $f$, denoted as $g_\theta(f)$ is the unique, strictly positive solution (if it exists) of the equation, e.g., for functions on $[0,S]$ or $\mathbb{R}^+$:

$$\frac{1}{x} \int_0^x f(s)\, ds = \theta x \tag{3}$$

again for $\theta > 0$. Also, the $g$-index is a very popular informetric tool in the impact evaluation of $f$, as is well-known.

Relations (2) and (3) show that they can be extended into the following defining equation for $x$. Let $A:(x,\theta) \rightarrow A(x,\theta)$ be an injective (in each variable separately), continuous, positive function on $\mathbb{R}^{+} \mathrm{x}\, \mathbb{R}_0^{+}$ (or subintervals of $\mathbb{R}^{+}$) and $T: \mathcal{V} \rightarrow \mathcal{W}$ an operator from a set of functions $\mathcal{V}$ to another set of functions $\mathcal{W}$ such that $\forall f \in \mathcal{V}$, $T(f)$ is positive and $T(f) = \mathbf{0} \Leftrightarrow f = \mathbf{0}$, where $\mathbf{0}$ is the null function. We further assume that for all $\theta$, the equation

$$T(f)(x) = A(x,\theta) \tag{4}$$

has a unique solution, denoted $x = m_\theta(f)$. Hence we define:

$$x = m_\theta(f) \Leftrightarrow T(f)(x) = A(x,\theta) \tag{5}$$

for all $\theta$ for which (5) exists (uniquely as assumed). We call this an admissible $\theta$ and they are defined by, for each $x$ in the domain of $A$ :

$$\theta =: \psi_f(x) = A^{-1}(T(f)(x)) \tag{6}$$

where $A^{-1}$ is (abuse of notation) the inverse function of $A(x,.): \theta \rightarrow A(x,\theta)$, which exists since $A(x,.)$ is injective.

The expression $m(f): \theta \rightarrow m_\theta(f)$ is a type of bundle in the sense of (Egghe & Rousseau, 2022), and, because of (1) and (2) we call the right-hand side of (5) a *Hirsch-type equation* but its solution is far more general, as the next examples show. The solution $m_\theta(f)$ we define as a *Hirsch-type bundle*.

Examples

(i) $T = \mathbf{1}$, the identical operator $T(f) = f$. By (5) we have

$$x = m_\theta(f) \Leftrightarrow f(x) = A(x,\theta) \tag{7}$$

They are the so-called “pointwise defined” bundles of the type (PED) as introduced in (Egghe, 2021). For $A(x,\theta) = \theta x$, we refind the $\theta$-dependent Hirsch function. For $A(x,\theta) = \theta x^p$ ($p \in \mathbb{R}^{+}$) we find the $\theta$-dependent Kosmulski indices (Kosmulski, 2006). Of course, for p =1, we have (2) again.

(ii) That the generalization in (5) with operator T is useful is seen as follows. Define for $x \geq a \geq 0$,

$$I(f)(x) = \int_a^x f(s)ds \tag{8}$$

where $f$ is defined on an interval starting at $a$. Then define $T = \mu$, for $x > a$:

$$\mu(f)(x) = \frac{1}{x-a} \int_a^x f(s)ds = \frac{I(f)(x)}{x-a} \tag{9}$$

and $\mu(f)(a) =: f(a)$ , since for continuous $f$ we have $\lim_{x \to a} \mu(f)(x) = f(a)$ by de l'Hôspital's rule. Now definition (5) becomes, for $T = \mu$:

$$x = m_\theta(f) \Leftrightarrow \mu(f)(x) = \frac{1}{x-a} \int_a^x f(s)ds = A(x, \theta) \tag{10}$$

Note that $\mu(f)$ decreases for decreasing $f$. For constant $f$ this reduces to (PED) bundles as in (7). For $A(x,\theta) = \theta(x-a)$ we have (3) (for $a = 0$), the $\theta$-dependent $g$-index $g_\theta(f)$. It is clear that, for $A(x,\theta) = \theta(x-a)^p, p > 0$, we find a $g$-index version of the Kosmulski indices (not yet defined as far as we are aware of). Note that (10) is equivalent with (for x > a)

$$I(f)(x) = (x-a)A(x,\theta) \tag{11}$$

which is the original definition of $g_\theta$. We prefer definition (10) though, since, for constant functions, we have (PED) as in (7) with the same $A(x,\theta)$ as in (10) and, more importantly, $\mu(f)$ decreases for decreasing $f$, while $I(f)$ increases.

In the next section, we provide important inequalities for the bundles of type (5), dependent on the behavior of $T$ and $A$. They yield the necessary tools for the study of two important informetric notions: convergence properties of the bundle $m$ in the following section and impact properties of $m$ in the last section before we come to the conclusions.

## 2. Basic inequalities for Hirsch-type equations (5)

We start with a lemma and two consequences for the bundle $m$ .

Lemma 1. For a function $k \in \mathcal{V}$ and for all admissible $\theta$ the following statements are true:

If $T(k)(.) - A(.,\theta)$ decreases, then

(i) $T(k)(x) > A(x,\theta) \Rightarrow m_\theta(k) > x$ and similarly: $T(k)(x) \geq A(x,\theta) \Rightarrow m_\theta(k) \geq x$

(ii) $T(k)(x) < A(x,\theta) \Rightarrow m_\theta(k) < x$ and similarly $T(k)(x) \leq A(x,\theta) \Rightarrow m_\theta(k) \leq x$

If $T(k)(.) - A(.,\theta)$ increases, then

(iii) $T(k)(x) > A(x,\theta) \Rightarrow m_\theta(k) < x$ and similarly: $T(k)(x) \geq A(x,\theta) \Rightarrow m_\theta(k) \leq x$

(iv) $T(k)(x) < A(x,\theta) \Rightarrow m_\theta(k) > x$ and similarly $T(k)(x) \leq A(x,\theta) \Rightarrow m_\theta(k) \geq x$

Hence all implications are equivalences.

Proof

(i) By (5) and since $m_\theta(k)$ exists and is unique, we have $x = m_\theta(k) \Leftrightarrow T(k)(x) = A(x,\theta)$. If $T(k)(x) > A(x,\theta)$ then we have for all $y < x$: $T(k)(y) - A(y,\theta) \geq T(k)(x) - A(x,\theta) > 0$ and hence $y \neq m_\theta(x)$. Hence $m_\theta(k) \geq x$. As, by (5), $m_\theta(k) \neq x$, we have $m_\theta(k) > x$. These two arguments together also yield $[T(k)(x) \geq A(x,\theta] \Rightarrow m_\theta(k) \geq x$ .

The proofs of the two other assertions are similar. □

Corollary 1. For functions $k,f \in \mathcal{V}$ we have:

(i) If $T(k)(.) - A(.,\theta)$ decreases and $T(k) > (\text{resp. } \geq) T(f)$ then $m_\theta(k) > (\text{resp. } \geq) m_\theta(f)$. A similar result holds for the case < (resp.) ≤.

(ii) If $T(k)(.) - A(.,\theta)$ increases and $T(k) >$ (resp. $\geq$) $T(f)$ then $m_\theta(k) <$ (resp. $\leq$) $m_\theta(f)$. A similar result holds for the case < (resp.) ≤.

Proof of (i)

We know that we have for all $x$ : $T(k)(x) >$ (resp.$\geq$) $T(f)(x)$. Then for $x = m_\theta(f)$ we have, by (5), $T(k)(x) >$ (resp. $\geq$) $T(f)(x) = A(x,\theta)$ and hence $m_\theta(k) >$ (resp. $\geq$) $m_\theta(f)$ by Lemma 1.

The proof of < (resp. ≤ ) in (i) and of (ii) is similar and left to the reader. □

Corollary 2. For a function $f \in \mathcal{V}$ we have:

If $T(k)(.) - A(.,\theta)$ decreases, then the functions $A(x,.)$ and $m_{(.)}(f)$ have opposite strict monotonicity:

(i) $A(x,.)$ strictly increasing implies $(\theta < \theta') \Rightarrow m_{\theta'}(f) < m_\theta(f)$

(ii) $A(x,.)$ strictly decreasing implies $(\theta < \theta') \Rightarrow m_{\theta'}(f) > m_\theta(f)$

If $T(k)(.) - A(.,\theta)$ increases, then the functions $A(x,.)$ and $m_{(.)}(f)$ have the same strict monotonicity:

(iii) $A(x,.)$ strictly increasing implies $(\theta < \theta') \Rightarrow m_\theta(f) < m_{\theta\prime}(f)$

(iv) $A(x,.)$ strictly decreasing implies $(\theta < \theta') \Rightarrow m_{\theta'}(f) < m_\theta(f)$

Proof (i)

$\theta < \theta'$ implies $A(x,\theta) < A(x,\theta')$. For $x = m_\theta(f)$ we have, by (5), $T(f)(x) = A(x,\theta)$. So $T(f)(x) < A(x,\theta')$. By (ii) in Lemma 1 (for $\theta'$): $m_{\theta'}(f) < m_\theta(f)$ .

For (ii) we use Lemma 1, part (i). For (iii) we use Lemma 1, part (iv) and for (iv) we use Lemma 1 part (iii). □

Note. For differentiable functions, Corollary 2 can be obtained by taking the derivative for $\theta$ in equation (5). We do not include this derivation here as our results are valid for non-differentiable functions as well.

Important comments on Lemma 1 and its consequences.

Here we take *T* = **1**, hence for all $k \in \mathcal{V}, T(k) = k$. We see from Corollary 1 that only in case (i) we have $[k > (\text{resp.} \geq) f] \Rightarrow [m_\theta(k) > (\text{resp.} \geq)\, m_\theta(f)]$, a necessary condition for an evaluation measure, such as for measuring impact, see further. Case (i) requires $k(.) - A(.,\theta)$ to be decreasing. For decreasing functions $k, f \in \mathcal{V}$ we illustrate this in terms of increasing $A(.,\theta)$, but also decreasing $A(.,\theta)$ is possible, see Figs 2(a),(b),(c).

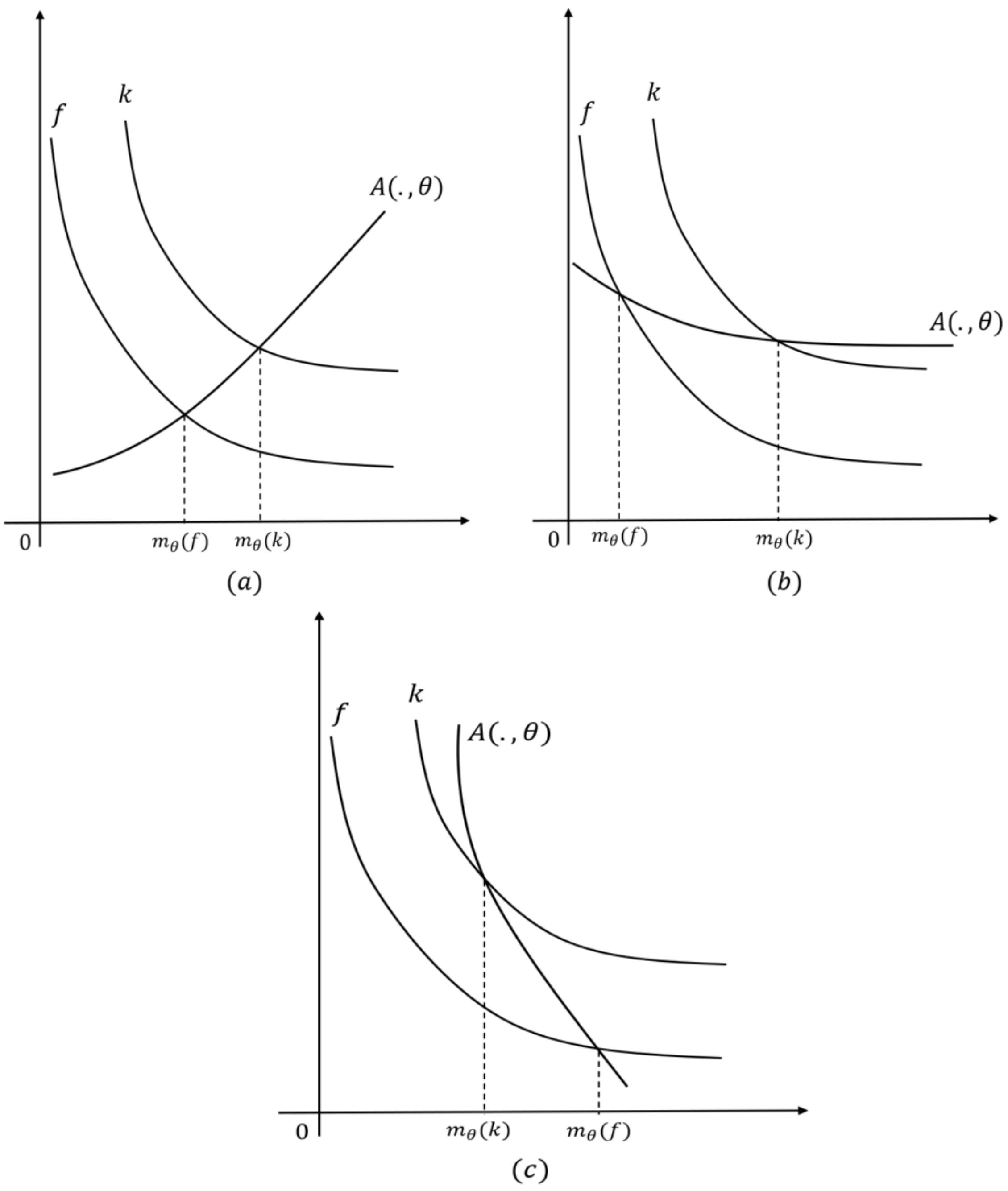

Fig.2. In cases (a) and (b), (i) in Corollary 1 is valid, while in case (c), (ii) is valid

Fig.2 (a) is the classical case, representing an increasing $A(.,\theta)$, hence (i) in Corollary 1. Examples: all (PED) cases, including Hirsch and Kosmulski cases, and when replacing $f$ and $k$ by $T(f)$ and $T(k)$ we also have all $g$-index cases. Also, a decreasing $A(.,\theta)$ is possible as long as $f$ and $k$ decrease faster than $A(.,\theta)$;

see Fig. 2(b). Also here we are in the (desired) case (i) in Corollary 1. Fig. 2(c) illustrates the (undesired) case (ii) in Corollary 1 where $A(.,\theta)$ decreases faster than $f$ and $k$ and hence $k(.) - A(.,\theta)$ increases. Similar figures can be drawn for increasing functions $f$ and $k$ but they are not studied here. In the same way, we can illustrate Corollary 2.

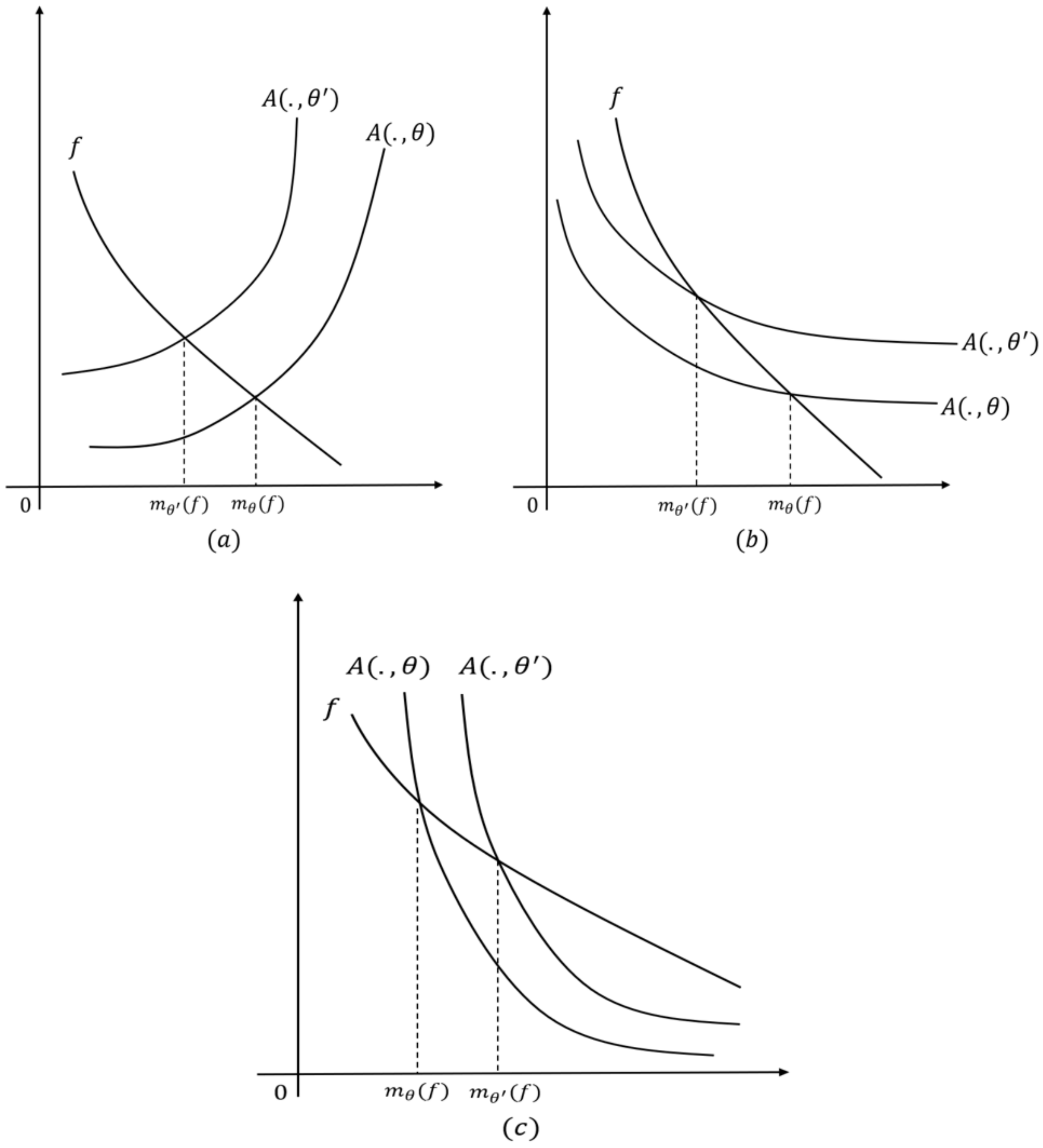

Fig. 3. In cases (a) and (b), (i) in Corollary 2 is valid and in case (c), (iii) is valid

The illustrations, Fig.3, are for $A(x,.)$ strictly increasing (but we draw $A(.,\theta)$, the same dependence as in $f$): $\theta < \theta' \Rightarrow A(.,\theta) < A(.,\theta')$. In (a) and (b), $f(.) - A(.,\theta)$ decreases: case (i) in Corollary (2) and in (c), $f(.) - A(.,\theta)$ increases: case (iii) in Corollary (2).

With these basic inequality results, we can now state and prove the Basic Theorem for Hirsch-type bundles (as well as equations).

Basic Theorem for Hirsch-type bundles

For all admissible $\theta$ and $(f_n)_{n\in\mathbb{N}}, f, (f_n, f \in \mathcal{V})$:

(i) If $T(f_n)(.)$ decreases for all $n \in \mathbb{N}$, and $A(.,\theta)$ increases, or, if $T(f_n)(.)$ increases for all $n \in \mathbb{N}$, and $A(.,\theta)$ decreases, then

$$|A(m_\theta(f_n),\theta) - A(m_\theta(f),\theta)| \le |T(f_n)(m_\theta(f) - T(f)(m_\theta(f)| \qquad (13)$$

(ii) If $T(f_n)(.)$ increases for all $n \in \mathbb{N}$, and $T(f_n)(.) - A(.,\theta)$ decreases, or if $T(f_n)(.)$ decreases for all $n \in \mathbb{N}$, and $T(f_n)(.) - A(.,\theta)$ increases, then

$$\left|T(f_n)\big(m_\theta(f)\big) - T(f)(m_\theta(f))\right| \le |A(m_\theta(f_n),\theta) - A(m_\theta(f),\theta)| \quad (14)$$

Proof. (i) First case. Since $A(.,\theta)$ increases we have: $|A(m_\theta(f_n),\theta) - A(m_\theta(f),\theta)|$

$$\begin{cases} = \text{(a), if } \mathrm{m}_\theta(f) \ge m_\theta(f_n), A(m_\theta(f),\theta) - A(m_\theta(f_n),\theta) \\ = (b), \text{if } \mathrm{m}_\theta(f) \le m_\theta(f_n), A(m_\theta(f_n),\theta) - A(m_\theta(f),\theta) \end{cases}$$

In case (a):

$0 \le (a) = T(f)\big(m_\theta(f)\big) - T(f_n)(m_\theta(f_n))$ (by (5))

$\le T(f)\big(m_\theta(f)\big) - T(f_n)(m_\theta(f))$ (since $T(f_n)(.)$ is decreasing ).

In case (b):

$0 \leq (b) = T(f_n)\big(m_\theta(f_n)\big) - T(f)(m_\theta(f))$ (by (5))

$\leq T(f_n)\big(m_\theta(f)\big) - T(f)(m_\theta(f))$ (since $T(f_n)(.)$ is decreasing ).

So, in both cases:

$|A(m_\theta(f_n),\theta) - A(m_\theta(f),\theta)| \leq |T(f_n)(m_\theta(f) - T(f)(m_\theta(f)|$ which proves (13) in the first case of (i). The proof of the second case of (i) is similar.

(ii) First case:

$\big|T(f_n)\big(m_\theta(f)\big) - T(f)(m_\theta(f))\big| = \big|T(f_n)\big(m_\theta(f)\big) - A(m_\theta(f),\theta)\big|$ , by (5)

$= (a), \quad if\ T(f_n)\big(m_\theta(f)\big) \geq A(m_\theta(f),\theta), \quad T(f_n)\big(m_\theta(f)\big) - A(m_\theta(f),\theta)$
$= (b), \quad if\ T(f_n)\big(m_\theta(f)\big) \leq A(m_\theta(f),\theta), \quad A(m_\theta(f),\theta) - T(f_n)(m_\theta(f))$

In case (a) we have, since $T(f_n)(.) - A(.,\theta)$ decreases it follows from the inequality in (a) and Lemma 1 (i) that $m_\theta(f_n) \geq m_\theta(f)$ and hence, since $T(f_n)(.)$ increases: $T(f_n)\big(m_\theta(f_n)\big) \geq T(f_n)\big(m_\theta(f)\big)$. But by (5), $T(f_n)(m_\theta(f_n)) = A(m_\theta(f_n),\theta)$. Hence

$$0 \leq (a) \leq A(m_\theta(f_n),\theta) - A(m_\theta(f),\theta)$$

In case (b):

Since $T(f_n)(.) - A(.,\theta)$ decreases, it follows from the inequality in (b) and Lemma 1 (ii) that $m_\theta(f_n) \leq m_\theta(f)$ and hence, since $T(f_n)(.)$ increases: $T(f_n)\big(m_\theta(f_n)\big) \leq T(f_n)\big(m_\theta(f)\big)$ . But, by (5), $T(f_n)(m_\theta(f_n)) = A(m_\theta(f_n),\theta)$ . Hence $0 \leq (b) \leq A(m_\theta(f),\theta) - A(m_\theta(f_n),\theta)$, whence (14)in both cases (a) and (b) which proves the first case of (ii). The proof of the second case is similar. □

Remark 1. (13) together with (14) imply equality in these formulae. However, the conditions in (i) and (ii) on $T$ and $A$ are mutually exclusive so we did not prove equality in (13) and (14). But inequalities as in (13) and (14) imply convergence results as will be proved in the next section. We can also say that (13) and (14) are stronger results than convergence and are also used in results on impact of the Hirsch-type bundles. This topic will be dealt with in the sequel.

Remark 2. In the Basic Theorem, we could have worked with only two functions $f$ and $g$, replacing the $f_n$, but since we want to use (13) and (14) in order to prove convergence results (in the next section) we formulated them in terms of $f$ and a sequence $(f_n)_{n\in\mathbb{N}}$.

## III. Convergence results for Hirsch-type bundles

Formulae (13) and (14) immediately imply the following convergence results:

Theorem 2. Assume that $\forall n \in \mathbb{N}$, $T(f_n)(.)$ is either increasing or decreasing, then:

(i) If $T$ is pointwise, sequentially continuous, and $A(.,\theta)$ has the opposite monotonicity of $T(f_n)(.)$, then, pointwise

$$(\lim_{n\to\infty} f_n = f) \Rightarrow (\lim_{n\to\infty} m(f_n) = m(f)) \qquad (15)$$

where $m(f)$ is the function $\theta \to m_\theta(f)$ and similarly for $m(f_n)$, for admissible $\theta$.

(ii) If $T^{-1}$ is pointwise, sequentially continuous and $T(f_n)(.) - A(.,\theta)$ has the opposite monotonicity of $T(f_n)(.)$, then, pointwise:

$(\lim_{n\to\infty} m(f_n) = m(f)) \Rightarrow (\lim_{n\to\infty} f_n = f)$ , where the equality in the consequent holds on the range of m(f) (16)

Proof

(i) We see that, because of the pointwise sequential continuity of $T$, the right-hand side of (13) goes to zero. Hence also the left-hand side: $\forall\theta$ admissible:

$$\lim_{n\to\infty} A(m_\theta(f_n),\theta) = A(m_\theta(f),\theta)$$

Since $A(.,\theta)$ is injective and continuous on (a subinterval of) $\mathbb{R}^+$, we have that its inverse exists and is also continuous. Hence (15) is proven.

(ii) Since $A(.,\theta)$ is continuous we have that the right-hand side of (14) goes to zero. Hence also the left-hand side of (14). By the pointwise sequential continuity of $T^{-1}$ (and since $m_\theta(f)$ is fixed) we have (16) on all $x$ of the form $x = m_\theta(f)$ for a certain admissible $\theta$. □

Corollary 1 (Egghe, 2022b).

If $(\lim_{n\to\infty} f_n = f)$ then

(i) $(\lim_{n\to\infty} m(f_n) = m(f))$ if m is (PED) (see (7)) with $A(.,\theta)$ increasing. This is hence also true for $m = h$.

(ii) $\lim_{n\to\infty} g(f_n) = g(f)$ with g as in (10)

Proof. (i) Here $T = \mathbf{1}$, $f_n(.)$ decreases $\forall n \in \mathbb{N}$. So (15) proves (i).

(ii) For $m = g, T = \mu, f_n(.)$ decreasing $\forall n \in \mathbb{N}$, it follows that $T(f_n)(.)$ is decreasing while $A(x,\theta) = x\,\theta$ (strictly increasing in $x$). Further, $T$ is pointwise sequentially continuous as follows from the theorem of Arzelà (for more on this, see (Egghe, 2022b)). Hence (15) proves (ii) □

Corollary 2 (see also (Egghe, 2022b))

For a (PED) such that $\lim_{n\to\infty} m(f_n) = m(f)$ pointwise, we have that $\lim_{n\to\infty} f_n = f$ pointwise if $f_n(.)$ decreases and $f_n(.) - A(.,\theta)$ increases.

Proof. Here again $T = \mathbf{1}$, and hence the result follows from (16). The result is also valid for $f_n(.)$ increasing and $f_n(.) - A(.,\theta)$ decreasing, but in this paper the functions $f_n$ and $f$ are supposed to be decreasing. □

Remark that both corollaries 1 and 2 are not only valid for the $h$- bundle (since it is a (PED)) but also for the polar form of a function $f$ since $\rho_\varphi(f) = h_\theta(f)\sqrt{1+\theta^2}$ with $\theta = tg(\varphi)$ fixed (see (Egghe, 2022b; Egghe & Rousseau, 2020).

In Egghe (2022b) also "uniform convergence" results are proved and also for this type of convergence, the inequalities (13) and (14) can be used. Indeed, if the right-hand side of (13) or (14) converges to zero uniformly, then, obviously, the left-hand side of (13) and (14) converges to zero uniformly. As shown in (Egghe, 2002) the uniform convergence of the $h$- and $g$-bundle (uniform in $\theta$) is obtained for $\theta > \theta_0 > 0$ ($\exists\, \theta_0$ fixed), i.e., for $\theta$ bounded away from zero. These results can be reproved using (13) and (14). We leave these exercises to the interested reader. But with (13) and (14) even more general results can be proved.

Remark. All these convergence results are important in informetrics. They guarantee that, when objects (e.g., authors) are close, say $f_n$ to $f$, then their evaluation measures $m(f_n)$ are also close to $m(f)$. This is necessary for evaluation tools to be valuable. One such evaluation is the measurement of impact. To this important informetric feature, we are devoting the next chapter. We will see that the main condition in Lemma 1 and its corollaries and in the Basic Theorem, the monotonicity of $T(k)(.) - A(.,\theta),\ k \in \mathcal{V}$ , plays a central role.

## IV. Impact in Hirsch-type bundles

For a general bundle $m: \theta \rightarrow m_\theta: f \rightarrow m_\theta(f)$ we determined that the following four properties define an impact bundle (Egghe & Rousseau, 2022). For functions $f, g \in \mathcal{V}$, defined on an interval [0, S] or $\mathbb{R}^+$:

AX.1: $f = 0 \Leftrightarrow m_\theta(f) = 0, \forall\theta$ admissible

AX.2: $f \leq g \Rightarrow m_\theta(f) \leq m_\theta(g), \forall\theta$ admissible

AX.3: $\forall a \in ]0, S[ : f <_a g \Rightarrow m(f) <_a m(g)$ where we define

$$f <_a g \Leftrightarrow f|_{[0,a]} < g|_{[0,a]} \quad (17)$$

and

$$m(f) <_a m(g) \Leftrightarrow m_\theta(f) < m_\theta(g), \forall\, \theta \in \psi_f([0,a]) \cup \psi_g([0,a]) \quad (18)$$

where $\psi_f$ and $\psi_g$ are the functions that relate $x = m_\theta(f)$ (resp. $m_\theta(g)$) with $\theta$ (the admissible values), in this article given by (6).

AX.4:
$\forall a \in ]0, S[: f = g \text{ on } [0,a] \Rightarrow \forall\theta \in \psi_f([0,a]) = \psi_g([0,a]): m_\theta(f) = m_\theta(g)$

In addition to the logical conditions on $T$ and $A$ in the introduction, we add the also logical conditions:

$\forall a \in ]0, S]$, $T_a$ is strictly increasing where we define

$$T_a\left(f|_{[0,a]}\right) = T(f)|_{[0,a]} \quad (19)$$

$f$and $T(f)$ are decreasing and $\forall\theta \neq 0: A(x,\theta) = 0 \Leftrightarrow x = 0$

We refer the reader to Egghe and Rousseau (2022) for arguments about why these axioms are determining impact. The most characteristic axiom is AX.3, showing the real nature of impact: it boils down to the following observation. The functions in $\mathcal{V}$ stand for rank-frequency functions, so with the highest production in the lower ranks $x$ , expressed by requiring $f$ to be decreasing. Now $f <_a g$ expresses that $g$ has higher production on $[0,a]$ and AX.3 requires that this is measured by the bundle $m$ so that, for the corresponding $\theta, m_\theta(f) < m_\theta(g)$.

We have the following impact theorem for Hirsch-type bundles.

Theorem 3. Given $f \in \mathcal{V}$, then

(i) $T(f)(.) - A(.,\theta)$decreases for all admissible $\theta$

$$\Rightarrow$$

(ii) *m* is an impact bundle

Proof.

AX.1 By the assumptions we have $\forall$ admissible $\theta, x = m_\theta(f) = 0 \Leftrightarrow T(f)(x) = 0$ , for all *x*, where we have assumed that all *x* in the domain of $T(f)$ yield an admissible $\theta$. Hence this is equivalent with $T(f) = 0 \Leftrightarrow f = 0$.

AX.2 $f \leq g \Rightarrow T(f) \leq T(g)$, by (19). By (i) we have by Corollary 1 (ii) that $m_\theta(f) \leq m_\theta(g)$, $\forall\theta$ admissible.

AX.3. $\forall a \in ]0,S[$, $f <_a g \Rightarrow f|_{[0,a]} < g|_{[0,a]}$ by (17) and hence, by (19), $\forall x \in [0,a]$; $T(f)(x) < T(g)(x)$. Suppose first that $\theta \in \psi_f([0,a])$ (part of the admissible $\theta$). Hence, by (5): for $x \in [0,a]$: $x = m_\theta(f) \Leftrightarrow T(f)(x) = A(x,\theta)$ and hence $T(g)(x) > A(x,\theta)$. By Lemma 1(i) we have $m_\theta(g) > x = m_\theta(f)$.

Similarly, for $\theta \in \psi_g([0,a])$ (hence admissible): by (5):

$$x = m_\theta(g) \Leftrightarrow T(g)(x) = A(x,\theta).$$

So now $T(f)(x) < A(x,\theta)$ and hence, by Lemma 1(ii): $m_\theta(f) < x = m_\theta(g)$ . Hence AX.3 is completely proven.

AX.4 $\forall a \in ]0,S[$: $f = g$ on $[0,a] \Rightarrow T(f) = T(g)$ on $[0,a]$ , by (19). By (6): $\psi_f|_{[0,a]} = \psi_g|_{[0,a]}$ and by (5): $m_\theta(f) = m_\theta(g)$ , $\forall\, \theta \in \psi_f([0,a]) = \psi_g([0,a])$ . □

Theorem 3 has also a reverse (ii) ⇒ (i) under an additional condition on $T$ which is satisfied for the examples as used in this paper: $T = \mathbf{1}, T = \mu$, ... We do not give the proof since it is quite intricate. It can be obtained from the author on demand.

So, under very reasonable conditions on $T$ and $A$, including all the examples, we have that a Hirsch-type bundle $m$ is an impact bundle iff $T(f)(.) - A(.,\theta)$ decreases for all admissible $\theta$. This (and the previous sections) show that the behavior of $T(f)(.) - A(.,\theta)$ is crucial for their use in informetrics and it shows again (cf. also (Egghe, 2022b; Egghe & Rousseau, 2022)) that the Hirsch-, g- and Kosmulski bundles are good informetric tools.

Corollary (Egghe and Rousseau, 2022)

(i) All Hirsch- and Kosmulski bundles are impact bundles

(ii) The g-bundle is an impact bundle. The same is true for its Kosmulski variant. (This is a new result).

Proof

(i) Here $T = \mathbf{1}$, hence $T(f) = f$ decreases. Furthermore $A(x,\theta) = \theta\, x^p, p > 0$ increases. Hence $T(f)(.) - A(.,\theta)$ decreases. So Theorem 3 gives that all Hirsch-bundles (*p=1*) and Kosmulski bundles (*p>0*) are impact bundles.

(ii) Here $T = \mu$ and hence $\mu(f){:}\, x \to \frac{1}{x}\int_0^x f$ is decreasing for f decreasing (well-known and easy to prove). Furthermore $A(x,\theta) = \theta x$ is increasing hence $T(f)(.) - A(.,\theta)$ decreases. So Theorem 3 gives that the g-bundle is an impact bundle. The same is true for the Kosmulski-type variant of $g$ : using $A(x,\theta) = \theta x^p, p > 0$.

Conclusions

The formalism

$$x = m_\theta(f) \Leftrightarrow T(f)(x) = A(x,\theta) \quad (5)$$

Defines a broad spectrum of bundles that can be used in the evaluation measurement (e.g., impact measurement) in informetrics, comprising all Hirsch-, Kosmulski- and g-type bundles.

It turns out that the property

$$\text{“}T(f)(.) - A(.,\theta) \text{ is monotone “} \quad (20)$$

Is fundamental in the study of the bundle (5) and especially the decreasing case in (20) guarantees the best informetric properties such as convergence and impact together, to which all the above-mentioned examples belong.

We invite the reader to further investigate (5) and (20) in relation with other good informetric properties such as the ones studied in (Egghe, 2022b) (several types of convergence), and Egghe & Rousseau (2023a,b) (global impact in relation with the non-normalized Lorenz curve and many variants). Lemma 1 is simple but crucial in all the results of this paper. It describes what happens if (5) is not satisfied and relates the sign of $T(f)(x) - A(x,\theta)$ to the relation between $x$ and $m_\theta(f)$ (the value for which $T(f)(.) - A(.,\theta) = 0$, by (5).

For the interested probabilist: inequalities (13) and (14) are somewhat similar to inequalities in probability theory between functions $f$ and their conditional expectations $E^\Sigma(f)$, where $\Sigma$ is a $\sigma$-algebra denoting how deep we are averaging the function $f$. In this comparison, $T(f)$ takes the place of such conditional expectation. Incidentally, for the $g-$index, $T$ is really an average as indicated above. We invite the interested probabilist to further investigate this similarity. They are different things but they have the same methodology in common: proving inequalities where the left-hand side and right-hand side deal with different mathematical functions. Then the convergence (to zero) of the left-hand side is guaranteed as soon as the right-hand side converges (to zero), see (Egghe, 1984, Chapter IV).

Acknowledgements. The author thanks Li Li (Beijing) for drawing excellent illustrations, and Ronald Rousseau for stimulating discussions.